\documentclass[10pt]{amsart}
\usepackage{amssymb}
\theoremstyle{plain}
\newtheorem{theorem}{Theorem}
\newtheorem{corollary}{Corollary}

\newtheorem{proposition}{Proposition}
\theoremstyle{example}
\newtheorem{example}{Example}
\theoremstyle{definition}

\theoremstyle{remark}

\numberwithin{equation}{section}

\setbox0=\hbox{$+$}
\newdimen\plusheight
\plusheight=\ht0
\def\+{\;\lower\plusheight\hbox{$+$}\;}

\setbox0=\hbox{$-$}
\newdimen\minusheight
\minusheight=\ht0
\def\-{\;\lower\minusheight\hbox{$-$}\;}

\setbox0=\hbox{$\cdots$}
\newdimen\cdotsheight
\cdotsheight=\plusheight
\def\cds{\lower\cdotsheight\hbox{$\cdots$}}

\begin{document}
\title[  Ramanujan,  continued fractions and  real numbers]
 {Ramanujan and the Regular Continued Fraction Expansion of Real Numbers }
\author{J. Mc Laughlin}
\address{Mathematics Department\\
 Trinity College\\
300 Summit Street, Hartford, CT 06106-3100}
\email{james.mclaughlin@trincoll.edu}
\author{ Nancy J. Wyshinski}
\address{Mathematics Department\\
       Trinity College\\
        300 Summit Street, Hartford, CT 06106-3100}
\email{nancy.wyshinski@trincoll.edu}

\date{July 4, 2003}
\begin{abstract}
In some recent papers, the authors considered regular continued fractions
of the form
\[
[a_{0};\underbrace{a,\cdots,a}_{m},   \underbrace{a^{2},\cdots,a^{2}}_{m},
 \underbrace{a^{3},\cdots,a^{3}}_{m}, \cdots ],
\]
where $a_{0} \geq 0$, $a \geq 2$ and $m \geq 1$ are integers.  The limits of
such continued fractions, for general $a$ and in the cases $m=1$ and $m=2$,
were given as ratios of certain infinite series.

However, these formulae can be derived from known facts about two
continued fractions of Ramanujan. Motivated by these observations,
we give alternative proofs of the results of the previous authors
for the cases $m=1$ and
$m=2$ and also use known results about other $q$-continued fractions
investigated by Ramanujan to
derive the limits of other infinite families of regular continued fractions.
\end{abstract}

\maketitle

\section{Introduction}
It is an interesting problem to try to find irrational numbers whose regular
continued fraction expansion contains predictable patterns and which can
be expressed in some other form.

The most familiar class of such numbers
comprises the  quadratic irrationalities, $ \alpha =p+q \sqrt{D}$,
where $p$ and $q$
are rational, $q \not = 0$ and $D$ is a non-square positive integer.
Such numbers have a regular continued fraction expansion which is ultimately
periodic:
\begin{equation*}
\alpha=[a_{0}; a_{1}, \cdots , a_{k}, \overline{b_{1}, \cdots , b_{n}}\,].
\end{equation*}

Another class consists of the \emph{Hurwitzian} continued fractions of the form
\begin{multline*}
[a_{0}; a_{1}, \cdots , a_{k},
 f_{1}(1), \cdots , f_{n}(1), f_{1}(2), \cdots , f_{n}(2),\cdots\,]\\
=:
[a_{0}; a_{1}, \cdots , a_{k},\overline{
 f_{1}(m), \cdots , f_{n}(m)}\,]_{m=1}^{\infty}.
\end{multline*}
Here the $f_{i}(x)$ are polynomials with rational coefficients taking only
positive integral values for integral $x\geq 1$ and at least one is non-constant.
The closed form for Hurwitzian continued fractions is not known in general.
This class contains numbers like
{\allowdisplaybreaks
\begin{multline*}
e^{2/(2n+1)} =
 [1; n, 12n+6, 5n+2,\\
 \overline{1,1,(6m+1)n+3m,\,(24m+12)n +12m+6,\,(6m+5)n+3m+2}\,]_{m=1}^{\infty}.
\end{multline*}
}
A third class is due to D.H. Lehmer \cite{L73}, who found closed forms for the
numbers represented by regular continued fractions whose partial quotients
were either terms in an arithmetic progression,
\[
[0;a,a+c,a+2c,a+3c,\cdots ],
\]
 or terms in two interlaced arithmetic
progressions,
\[
[0;a,b,a+c,b+d,a+2c,b+2d,\cdots].
\]
 An example that Lehmer gave of the former type was the following:
\begin{equation*}
[1;2,3,4,5,\cdots]= \frac{\sum_{m=0}^{\infty}\frac{1}{(m!)^{2}}}
{\sum_{m=0}^{\infty}\frac{1}{m!(m+1)!}}
\end{equation*}
Tasoev \cite{T84}, \cite{T00} proposed a new type of continued fraction
of the form
\begin{equation}\label{taseq}
[a_{0};\underbrace{a,\cdots,a}_{m},   \underbrace{a^{2},\cdots,a^{2}}_{m},
 \underbrace{a^{3},\cdots,a^{3}}_{m}, \cdots ],
\end{equation}
where $a_{0} \geq 0$, $a \geq 2$ and $m \geq 1$ are integers.
This type  was further investigated by Komatsu  \cite{K03a},
who derived a closed form for the general case
($m\geq 1$, arbitrary).  For the special cases $m=1$ and $m=2$, he derived
the following expressions. As usual, the empty product denotes $1$.
\begin{equation}\label{k1}
[0,\overline{a^{k}}]_{k=1}^{\infty}:=
[0;a,a^{2},a^{3},a^{4},\cdots]=
\frac{\sum_{s=0}^{\infty}a^{-(s+1)^{2}}\prod_{i=1}^{s}(a^{2i}-1)^{-1}}
{\sum_{s=0}^{\infty}a^{-s^{2}}\prod_{i=1}^{s}(a^{2i}-1)^{-1}}.
\end{equation}
\begin{equation}\label{k2}
[0,\overline{a^{k},a^{k}}]_{k=1}^{\infty}:=
[0;a,a,a^{2},a^{2},\cdots]=
\frac{\sum_{s=0}^{\infty}a^{-(s+1)(s+2)/2}\prod_{i=1}^{s}(a^{i}-1)^{-1}}
{\sum_{s=0}^{\infty}a^{-s(s+1)/2}\prod_{i=1}^{s}(a^{i}-1)^{-1}}.
\end{equation}
In \cite{K03b} he generalized these results and gave similar
expressions for continued fractions like $[0,\overline{ua^{k}}]$
and $[0,ua-1,\overline{1,ua^{k+1}-2}]$, with $a>1$ an integer and
$u$ rational such that $ua \in \mathbb{Z}^{+}$.

Komatsu manipulated certain infinite series to derive his results.
 In this present paper we use known facts
about certain $q$-continued fractions studied by Ramanujan to
give alternative, perhaps simpler,  derivations of some of
Komatsu's results.

Motivated by this connection between
families of regular continued fraction expansions and
$q$-continued fractions,  we investigated other $q$-continued
fractions studied by Ramanujan and were able to derive other infinite families of
regular continued fraction expansions which can be summed in a
closed form. Here are some examples of our results (proofs are found
throughout the paper):
\begin{example}
If $c>1$, $a>1$ are  integers and $d$ is rational such that $da>1$, then
\begin{equation*}
[c-1;\overline{1,d a^{k}-1, c a^{k}-1}\,]_{k=1}^{\infty}=
\frac{\sum_{n=0}^{\infty}
\displaystyle{
\frac{(-1)^{n}}
{c^{n-1}d^{n}a^{n^{2}}(-1/a;-1/a)_{n}}}
}
      {\sum_{n=0}^{\infty}
\displaystyle{
\frac{1}
{(cd)^{n}a^{n^{2}+n}(-1/a;-1/a)_{n}}
}}.
\end{equation*}
\end{example}
\begin{example}
Let $a>1$ be an integer.
Then
{\allowdisplaybreaks
\begin{equation*}
[1;\overline{a^{2k-1}-1,a^{2k}+1}]_{k=1}^{\infty}=
\frac{(-1/a;-1/a^{3})_{\infty}}
{(1/a^{2};-1/a^{3})_{\infty}}.
\end{equation*}
}
\end{example}
\begin{example}
Define
\[
F(c,d,q):=\sum_{n=0}^{\infty}\frac{(-1)^{n}c^{n}q^{n(n+1)/2}}{(q;q)_{n}(cq/d;q)_{n}}
\]
and let $\omega = e^{2\pi i/3}$. If $a>1$ is an integer and $c$ is a rational
such that $a/c$ is an integer, $a/c>2$, then
\begin{equation*}
\left [0;\frac{a}{c}-2,\overline{1,\frac{a^{k+1}}{c}-3}\right ]_{k=1}^{\infty}=
\frac{c/a\,F(-c \omega /a, \omega^{2},1/a)}
{(1+c \omega^{2}/a)F(-c \omega , \omega^{2},1/a)}.
\end{equation*}
\end{example}
\begin{example}
For $r$, $s$ and $q \in \mathbb{C}$ with $|q|<1$, define
\[
\phi(r,s,q)= \sum_{n=0}^{\infty}\frac{q^{(n^{2}+n)/2}r^{n}}
{(q;q)_{n}(-sq;q)_{n}}.
\]
Let $m$ and $n$ be positive integers and let $d$ be rational such that
$d n \in \mathbb{Z^{+}}$ and $d m n>1$.
If $n>2$ and $m>1$ then
\begin{multline*}
[0,\overline{1,d^{2k-2}n^{2k-1}-2,1,m^{2k-1}-1,
d^{2k-1}n^{2k},m^{2k}-1
}]_{k=1}^{\infty}\\
=1+
\frac{\phi(d m,d,-1/(d m n))}
{\phi(-1/n,d,-1/(d m n))}.
\end{multline*}
\end{example}

Curiously, it seems that Ramanujan was not particularly interested in
the regular continued fraction expansion of real numbers, or with his flair
for continued fractions he would almost certainly
have derived our results himself.
Perhaps he typified a gap that may exist between people who study
$q$-continued fractions and those who are interested
in the regular continued fraction expansion of real numbers?

\section{Some Continued Fractions of Ramanujan}
 In Chapter 16 of the Second Notebook, Ramanujan
gave the following corollary to Entry 15 (\cite{B91}, page 30):
\begin{corollary}\label{ramcor}
If $|q|<1$, then
\begin{equation}\label{rrgeneq}
 G(b,q):= 1+
 \frac{bq}{1}
\+
 \frac{bq^{2}}{1}
\+
 \frac{bq^{3}}{1}
\+\,\cds=\frac{\sum_{n=0}^{\infty}
\displaystyle{
\frac{b^{n}q^{n^{2}}}{(q;q)_{n}}}
}
      {\sum_{n=0}^{\infty}
\displaystyle{
\frac{b^{n}q^{n^{2}+n}}{(q;q)_{n}}
}
}.
\end{equation}
\end{corollary}
We have restated this corollary in a form that is more convenient for our
purposes.
If we set $b=1$ then the well-known expression for the Rogers-Ramanujan
continued fraction follows.
We also recall the
famous Rogers-Ramanujan identities \cite{R94}:
\begin{align*}
\sum_{n=0}^{\infty}\frac{q^{n^{2}}}{(q;q)_{n}}&=
\prod_{j=0}^{\infty}\frac{1}{(1-q^{5j+1})(1-q^{5j+4})}=
\frac{1}{(q;q^{5})_{\infty} (q^{4};q^{5})_{\infty}},\\
&\phantom{as}\\
\sum_{n=0}^{\infty}\frac{q^{n^{2}+n}}{(q;q)_{n}}&=
\prod_{j=0}^{\infty}\frac{1}{(1-q^{5j+2})(1-q^{5j+3})}=
\frac{1}{(q^{2};q^{5})_{\infty} (q^{3};q^{5})_{\infty}},
\end{align*}
where
\[
(c;q)_{0}:=1,\,\,\,\, (c;q)_{n}=\prod_{j=0}^{n-1}(1-cq^{j}), \,\,\,
(c;q)_{\infty}=\lim_{n \to \infty}(c;q)_{n},\,\,\,
\text{  } |q|<1.
\]
If we combine these identities we have that
\begin{equation}\label{rreq}
K(q):=
1+
 \frac{q}{1}
\+
 \frac{q^{2}}{1}
\+
 \frac{q^{3}}{1}
\+\,\cds=
\frac{\sum_{n=0}^{\infty}
\displaystyle{
\frac{q^{n^{2}}}{(q;q)_{n}}}
}
      {\sum_{n=0}^{\infty}
\displaystyle{
\frac{q^{n^{2}+n}}{(q;q)_{n}}}
}
=\frac{(q^{2};q^{5})_{\infty} (q^{3};q^{5})_{\infty} }
{(q;q^{5})_{\infty} (q^{4};q^{5})_{\infty}}.
\end{equation}
If we transform $K(q)$ so that
\begin{equation*}
K(q):=
1+
 \frac{1}{1/q}
\+
\frac{1}{1/q}
\+
 \frac{1}{1/q^{2}}
\+
\frac{1}{1/q^{2}}
\+
 \frac{1}{1/q^{3}}
\+
 \frac{1}{1/q^{3}}
\+\,\cds
\end{equation*}
and set $q=1/a$, where $a \geq 2$ is a positive integer,
we have that
\begin{align*}
1+[0;a,a,a^{2},a^{2},\cdots]&=
\frac{\sum_{n=0}^{\infty}
\displaystyle{
\frac{1}{a^{n^{2}}(1/a;1/a)_{n}}}
}
      {\sum_{n=0}^{\infty}
\displaystyle{
\frac{1}{a^{n^{2}+n}(1/a;1/a)_{n}}}
}\\
&=\frac{(1/a^{2};1/a^{5})_{\infty} (1/a^{3};1/a^{5})_{\infty} }
{(1/a;1/a^{5})_{\infty} (1/a^{4};1/a^{5})_{\infty}}.
\end{align*}
The first equation, after some elementary manipulation, gives the identity
at \eqref{k2}.

We will use the identity at \eqref{rrgeneq} to prove the following theorem.

\begin{theorem}\label{t1}
Let $a\geq 2$ be a positive integer and suppose $c$ and $d$ are rationals such that
$ca$, $da \in \mathbb{Z^{+}}$. Then
\begin{equation}\label{eqgen1}
[0;\overline{d a^{k}, c a^{k}}\,]_{k=1}^{\infty}=
\frac{\sum_{n=0}^{\infty}
\displaystyle{
\frac{1}
{c^{n}d^{n+1}a^{(n+1)^{2}}(1/a;1/a)_{n}}}
}
      {\sum_{n=0}^{\infty}
\displaystyle{
\frac{1}
{(cd)^{n}a^{n^{2}+n}(1/a;1/a)_{n}}
}
}.
\end{equation}
\begin{equation}\label{eqgen2}
[0;\overline{c a^{k}}\,]_{k=1}^{\infty}=
\frac{\sum_{n=0}^{\infty}
\displaystyle{
\frac{1}
{c^{2n+1}a^{2n^2+3n+1}(1/a^{2};1/a^{2})_{n}}}
}
      {\sum_{n=0}^{\infty}
\displaystyle{
\frac{1}
{c^{2n}a^{2n^{2}+n}(1/a^{2};1/a^{2})_{n}}
}
}.
\end{equation}
If $ca>1$, then
\begin{equation}\label{geneq3}
[0;ca-1,\overline{1,ca^{k+1}-2}\,]_{k=1}^{\infty}=
\frac{\sum_{n=0}^{\infty}
\displaystyle{
\frac{(-1)^{n}}
{c^{2n+1}a^{2n^2+3n+1}(1/a^{2};1/a^{2})_{n}}}
}
      {\sum_{n=0}^{\infty}
\displaystyle{
\frac{(-1)^{n}}
{c^{2n}a^{2n^{2}+n}(1/a^{2};1/a^{2})_{n}}
}
}.
\end{equation}
If $da>1$ and $ca>2$, then
\begin{multline}\label{geneq4}
[0;da-1,1,ca-2,\overline{1,da^{k+1}-2,1,ca^{k+1}-2}\,]_{k=1}^{\infty}\\
=
\frac{\sum_{n=0}^{\infty}
\displaystyle{
\frac{(-1)^{n}}
{c^{n}d^{n+1}a^{(n+1)^{2}}(1/a;1/a)_{n}}}
}
      {\sum_{n=0}^{\infty}
\displaystyle{
\frac{(-1)^{n}}
{(cd)^{n}a^{n^{2}+n}(1/a;1/a)_{n}}
}
}.
\end{multline}
If $c>1$ is an integer and $da>1$, then
\begin{equation}\label{eqgen5}
[c-1;\overline{1,d a^{k}-1, c a^{k}-1}\,]_{k=1}^{\infty}=
\frac{\sum_{n=0}^{\infty}
\displaystyle{
\frac{(-1)^{n}}
{c^{n-1}d^{n}a^{n^{2}}(-1/a;-1/a)_{n}}}
}
      {\sum_{n=0}^{\infty}
\displaystyle{
\frac{1}
{(cd)^{n}a^{n^{2}+n}(-1/a;-1/a)_{n}}
}
}.
\end{equation}
\end{theorem}
\begin{proof}
On the left side of Equation \ref{rrgeneq}, let $b=1/(c\,d)$, so that $G(b,q)$ can be written
as
\begin{equation*}
G(b,q)=1 +
\frac{1}{c}
\left (
\frac{1}{d/q}
\+
\frac{1}{c/q}
\+
\frac{1}{d/q^{2}}
\+
\frac{1}{c/q^{2}}
\+
\cds
\right ).
\end{equation*}
If we now let $q=1/a$ and $c$ and $d$ be rationals such that
$c a$, $d a \in \mathbb{Z^{+}}$
then
\begin{equation}\label{eqgen0}
1+\frac{1}{c}[0;\overline{d a^{k}, c a^{k}}\,]_{k=1}^{\infty}=
\frac{\sum_{n=0}^{\infty}
\displaystyle{
\frac{1}
{(cd)^{n}a^{n^{2}}(1/a;1/a)_{n}}}
}
      {\sum_{n=0}^{\infty}
\displaystyle{
\frac{1}
{(cd)^{n}a^{n^{2}+n}(1/a;1/a)_{n}}
}
}.
\end{equation}
After some manipulation, we have the identity at \eqref{eqgen1}.

Remark: Komatsu also has this result in \cite{K03b} as a corollary to a more general
result about continued fractions of the form
$[0;\overline{d a^{k}, c b^{k}}\,]_{k=1}^{\infty}$. If we set
$c=d=1$, some further simplification
gives the result from his paper \cite{K03a} at  \eqref{k2}.

If we replace $a$ by $a^{2}$ and $d$ by $c/a$ in Equation \ref{eqgen1}, we have the
identity at \eqref{eqgen2}.
This is Theorem 1 from \cite{K03b} and setting $c=1$ gives the identity
at \ref{k1} from \cite{K03a}.

As shown in \cite{VDP94}, negatives and zeroes can easily be removed from
regular continued fraction expansions. Indeed, it is easy to check that
$[m,n,0,p, \alpha]$ $=$ $[m,n+p, \alpha]$ and $[m,-n,\alpha]=[m-1,1,n-1,-\alpha]$.
Hence it is possible to allow the parameters $a$, $c$ and $d$ in Equations
\ref{eqgen0},  \ref{eqgen1} and \ref{eqgen2} to take negative values. Thus, if we
replace $c$ by $-c$ and $a$ by $-a$ in \eqref{eqgen2}, for example,
and repeatedly apply the second
of the above conditions, we have that
\begin{align*}
[0&;ca,-ca^2,ca^3,-ca^4,ca^5,-ca^6,\cdots]\\
&=[0;ca-1,1,ca^2-1,-ca^3,ca^4,-ca^5,ca^6,\cdots]\\
&=[0;ca-1,1,ca^2-2,1,ca^3-1,-ca^4,ca^5,-ca^6,\cdots]\\
&=[0;ca-1,1,ca^2-2,1,ca^3-2,1,ca^4-1,-ca^5,ca^6,\cdots]\\
&=[0;ca-1,1,ca^2-2,1,ca^3-2,1,ca^4-2,1,ca^5-1,-ca^6,\cdots]
\text{ etc.}
\end{align*}
Finally, we have, for $ca>1$, the identity at \eqref{geneq3}

This is Theorem 2 from \cite{K03b}. Likewise, if $c$ is replaced by $-c$ in \ref{eqgen1}
and the resulting continued fraction similarly manipulated, one gets \eqref{geneq4}.
This is Theorem 4 from \cite{K03b}.
One that Komatsu missed is \eqref{eqgen5}.
This is derived from Equation \ref{eqgen0} by multiplying
across by $c$, replacing $a$ by $-a$ and similarly manipulating the continued fraction
to remove the negatives.
\end{proof}
We digress slightly before proving our next result. We introduce some notation
from \cite{LW92}
(page 83).
We call $d_{0}+K_{n=1}^{\infty}c_{n}/d_{n}$  a \emph{canonical contraction} of
 $b_{0}+K_{n=1}^{\infty}a_{n}/b_{n}$ if
\begin{align*}
&C_{k}=A_{n_{k}},& &D_{k}=B_{n_{k}}& &\text{ for } k=0,1,2,3,\ldots \, ,\phantom{asdasd}&
\end{align*}
where $C_{n}$, $D_{n}$, $A_{n}$ and $B_{n}$ are canonical numerators and denominators
of $d_{0}+K_{n=1}^{\infty}c_{n}/d_{n}$ and $b_{0}+K_{n=1}^{\infty}a_{n}/b_{n}$ respectively.

From \cite{LW92} (page 85)we also have the following theorem.
\begin{theorem}\label{odcf}
The canonical contraction of $b_{0}+K_{n=1}^{\infty}a_{n}/b_{n}$ with
$C_{0}=A_{1}/B_{1}$
\begin{align*}
&C_{k}=A_{2k+1}& &D_{k}=B_{2k+1}& &\text{ for } k=1,2,3,\ldots \, ,&
\end{align*}
exists if and only if $b_{2k+1} \not = 0 for K=0,1,2,3,\ldots$, and in this case is given by
{\allowdisplaybreaks
\begin{multline}\label{E:odcf}
\frac{b_{0}b_{1}+a_{1}}{b_{1}}
-
\frac{a_{1}a_{2}b_{3}/b_{1}}{b_{1}(a_{3}+b_{2}b_{3})+a_{2}b_{3}}
\-
\frac{a_{3}a_{4}b_{5}b_{1}/b_{3}}{a_{5}+b_{4}b_{5}+a_{4}b_{5}/b_{3}}\\
\-
\frac{a_{5}a_{6}b_{7}/b_{5}}{a_{7}+b_{6}b_{7}+a_{6}b_{7}/b_{5}}
\-
\frac{a_{7}a_{8}b_{9}/b_{7}}{a_{9}+b_{8}b_{9}+a_{8}b_{9}/b_{7}}
\+
\cds .
\end{multline}
}
\end{theorem}
The continued fraction \eqref{E:odcf} is called the \emph{odd} part of $b_{0}+K_{n=1}^{\infty}a_{n}/b_{n}$.
The following corollary follows easily from Theorem \ref{odcf}.
{\allowdisplaybreaks
\begin{corollary}\label{corcf7}
The odd part of the continued fraction
%{\allowdisplaybreaks
\begin{equation*}
\frac{c_{1}}{1}
\-
\frac{c_{2}}{1}
\+
\frac{c_{2}}{1}
\-
\frac{c_{3}}{1}
\+
\frac{c_{3}}{1}
\-
\frac{c_{4}}{1}
\+
\frac{c_{4}}{1}
\-
\cds
\end{equation*}
%}
is
%{\allowdisplaybreaks
\begin{equation*}
c_{1}
+
\frac{c_{1}c_{2}}{1}
\+
\frac{c_{2}c_{3}}{1}
\+
\frac{c_{3}c_{4}}{1}
\+
\cds .
\end{equation*}
%}
\end{corollary}
}

We will also make use of Worpitzky's Theorem  (see \cite{LW92}, pp. 35--36).
\begin{theorem}(Worpitzky)
 Let the continued fraction $K_{n=1}^{\infty}a_{n}/1$ be such that
$|a_{n}|\leq 1/4$ for $n \geq 1$. Then$K_{n=1}^{\infty}a_{n}/1$ converges.
 All approximants of the continued fraction lie in the disc $|w|<1/2$ and the value of the
continued fraction is in the disk $|w|\leq1/2$.
\end{theorem}
We are now ready to prove Theorem \ref{t2}

\begin{theorem}\label{t2}
Let $m$, $r$, $s$ and $a>1$ be positive integers such that $sa/(rm^{2})$
is an integer. If $m>1$ and $sa/(rm^{2})>1$, then
{\allowdisplaybreaks
\begin{multline}\label{eqex2}
\left[0;\overline{m-1,1,\frac{sa^{k}}{rm^2}-1,m-1,1,a^{k}-1}\,\right]_{k=1}^{\infty}\\
=\frac{1}{m}+\frac{\sum_{n=0}^{\infty}
\displaystyle{
\frac{(r/s)^{n+1}}{a^{(n+1)^{2}}(1/a;1/a)_{n}}}
}
      {\sum_{n=0}^{\infty}
\displaystyle{
\frac{(r/s)^{n}}{a^{n^{2}+n}(1/a;1/a)_{n}}
}
}.
\end{multline}
 If $m>2$ and $sa/(rm^{2})>1$, then
\begin{equation}\label{eqex3}
\left[1;\overline{m,\frac{sa^{k}}{rm^2}-1,1,m-2,1,a^{k}-1}\,\right]_{k=1}^{\infty}
=\frac{1}{m}+\frac{\sum_{n=0}^{\infty}
\displaystyle{
\frac{(-r/s)^{n}}{a^{n^{2}}(1/a;1/a)_{n}}}
}
      {\sum_{n=0}^{\infty}
\displaystyle{
\frac{(-r/s)^{n}}{a^{n^{2}+n}(1/a;1/a)_{n}}
}
}
.
\end{equation}
 If $m>2$ and $sa/(rm^{2})>1$, then
\begin{equation}\label{eqex4}
\left[0;\overline{1,m-2,1,\frac{sa^{k}}{rm^2}-1,m,a^{k}-1}\,\right]_{k=1}^{\infty}
=\frac{-1}{m}+\frac{\sum_{n=0}^{\infty}
\displaystyle{
\frac{(-r/s)^{n}}{a^{n^{2}}(1/a;1/a)_{n}}}
}
      {\sum_{n=0}^{\infty}
\displaystyle{
\frac{(-r/s)^{n}}{a^{n^{2}+n}(1/a;1/a)_{n}}
}
}
.
\end{equation}
 If $m>2$ and $sa/(rm^{2})>1$, then
\begin{multline}\label{eqex5}
\phantom{as}\\
\text{\small{$
\left[1;\overline{m,\frac{sa^{2k-1}}{rm^2}-1,1,m-1,a^{2k-1}-1,1,
m-2,1,\frac{sa^{2k}}{rm^2}-1,m-1,1,a^{2k}-1}\,\right]_{k=1}^{\infty}$}}\\
=\frac{1}{m}+
\frac{\sum_{n=0}^{\infty}
\displaystyle{
\frac{(-r/s)^{n}}{a^{n^{2}}(-1/a;-1/a)_{n}}}
}
      {\sum_{n=0}^{\infty}
\displaystyle{
\frac{(r/s)^{n}}{a^{n^{2}+n}(-1/a;-1/a)_{n}}
}
}
.
\end{multline}
}
\end{theorem}

\begin{proof}
Let $c_{1}$ be arbitrary and apply  Corollary \ref{corcf7}
to $G(b,q)$ from Corollary
\ref{ramcor}, so that
{\allowdisplaybreaks
\begin{align}\label{ramcor2}
G(b,q)=
1-c_{1}+
\frac{c_{1}}{1}
&\-
\frac{b q/c_{1}}{1}
\+
\frac{b q/c_{1}}{1}
\-
\frac{qc_{1}}{1}
\+
\frac{qc_{1}}{1}\\
&\-
\frac{b q^{2}/c_{1}}{1}
\+
\frac{b q^{2}/c_{1}}{1}
\-
\frac{q^{2}c_{1}}{1}
\+
\frac{q^{2}c_{1}}{1}
\-
\cds  \notag
\\
=
1-c_{1}+
\frac{1}{1/c_{1}}
&\-
\frac{1}{c_{1}^2/(b q)}
\+
\frac{1}{1/c_{1}}
\-
\frac{1}{1/q}
\+
\frac{1}{1/c_{1}} \notag
\\
&\-
\frac{1}{c_{1}^2/(b q^{2})}
\+
\frac{1}{1/c_{1}}
\-
\frac{1}{1/q^{2}}
\+
\frac{1}{1/c_{1}}
\-\notag
 \\
&\phantom{sadasdssada}\cds  \notag
 \\
&\-
\frac{1}{c_{1}^2/(b q^{n})}
\+
\frac{1}{1/c_{1}}
\-
\frac{1}{1/q^{n}}
\+
\frac{1}{1/c_{1}}
\-
\cds .\notag\\
\end{align}
}

The first equality is valid since the continued fraction on the left converges by
Worpitzky's Theorem
 and hence equals its odd part, which is $G(b,q)$, by
Corollary \ref{corcf7}.

If we now let $c_{1}=1/m$, where $m$ is an  integer,
$b=r/s$, where $r$ and $s$ are integers, and $q=1/a$, where $a$ is an integer
such that $sa/(r m^2)$ is an integer, we have, after a little manipulation to bring all
the negative signs into the denominators, that
{\allowdisplaybreaks
\begin{align}\label{eqex}
G(r/s,1/a)&=
1-\frac{1}{m} +[0; \overline{m,-s a^{k}/(r m^2),-m, a^{k}}\,]_{k=1}^{\infty}\\
&=\frac{\sum_{n=0}^{\infty}
\displaystyle{
\frac{(r/s)^{n}}{a^{n^{2}}(1/a;1/a)_{n}}}
}
      {\sum_{n=0}^{\infty}
\displaystyle{
\frac{(r/s)^{n}}{a^{n^{2}+n}(1/a;1/a)_{n}}
}
}.\notag
\end{align}
}

One can now make various choices for the signs of the parameters $m$, $s$, $r$ and
$a$ and remove the negative signs as before to produce regular continued fraction
expansions.
If we choose all parameters to be positive, then
 clearing negatives and a little manipulation on the right side of Equation \ref{eqex}
gives \eqref{eqex2}.

Remark:  We point out something that is a little curious.
Suppose we choose $r=s=1$ in \eqref{eqex2} and
suppose further that $m_{1}^2|a$ and $m_{2}^2|a$.
Then
{\allowdisplaybreaks
\begin{multline}
\left[0;\overline{m_{1}-1,1,\frac{a^{k}}{m_{1}^2},m_{1}-1,1,a^{k}-1}\,\right]_{k=1}^{\infty}\\
=\left[0;\overline{m_{2}-1,1,\frac{a^{k}}{m_{2}^2},m_{2}-1,1,a^{k}-1}\,\right]_{k=1}^{\infty}
+\frac{1}{m_{1}}-\frac{1}{m_{2}}.
\end{multline}
}
In other words, different values of the parameter $m$ produces
continued fractions whose values differ by a rational number. This
is clear since the only explicit appearance of $m$ on the right of
\eqref{eqex2}
 is in the fraction $1/m$.

If we replace $r$ by  $-r$  in \eqref{eqex}, then after clearing negatives we get
\eqref{eqex3}.\\
Remark:  If we let $m=2$ in \eqref{eqex3}, then the resulting zeroes can be removed
as described in the proof of Theorem \ref{t1}, to give that, if $4r|sa$ and $sa/4r>1$, then
 \begin{equation}\label{eqex3a}
\left[1;\overline{2,\frac{sa^{k}}{4r}-1,2,a^{k}-1}\,\right]_{k=1}^{\infty}
=\frac{1}{2}+\frac{\sum_{n=0}^{\infty}
\displaystyle{
\frac{(-r/s)^{n}}{a^{n^{2}}(1/a;1/a)_{n}}}
}
      {\sum_{n=0}^{\infty}
\displaystyle{
\frac{(-r/s)^{n}}{a^{n^{2}+n}(1/a;1/a)_{n}}
}
}
.
\end{equation}
Similar transformations of some of the other continued fractions in the paper are
possible. However, we generally ignore this.

If we replace  $m$ by $-m$ and  $s$ by $-s$ in \eqref{eqex},
then clearing negatives gives
\eqref{eqex4}. Finally, replacing $a$ by $-a$ \eqref{eqex} and clearing negatives gives
\eqref{eqex5}.

Other choices of signs will lead to different regular continued fraction expansions
at \eqref{eqex}.
\end{proof}

On page 290 of his second notebook, Ramanujan recorded the following
continued fraction identity.
\begin{proposition}$($see \cite{B98}, page 46, Entry 19$)$
For  $|q|<1$,
\begin{equation}\label{r2eq}
\frac{(q^{2};q^{3})_{\infty}}
{(q;q^{3})_{\infty}}
=
\frac{1}{1}
\-
\frac{q}{1+q}
\-
\frac{q^{3}}{1+q^{2}}
\-
\frac{q^{5}}{1+q^{3}}
\-
\cds .
\end{equation}
\end{proposition}
A proof of this identity can be found in \cite{ABSYZ03}.
We will use \eqref{r2eq} to prove the following theorem.
\begin{theorem}\label{t3}
Let $a>1$ be an integer.
Then
{\allowdisplaybreaks
\begin{equation}\label{r2a}
[1;\overline{a^{2k-1}-1,a^{2k}+1}]_{k=1}^{\infty}=
\frac{(-1/a;-1/a^{3})_{\infty}}
{(1/a^{2};-1/a^{3})_{\infty}}.
\end{equation}
}
{\allowdisplaybreaks
\begin{equation}\label{r2b}
[1;\overline{a^{k}-1,1}]_{k=1}^{\infty}=
\frac{(1/a^{2};1/a^{3})_{\infty}}
{(1/a;1/a^{3})_{\infty}}.
\end{equation}
}
\end{theorem}
\begin{proof}
Replace $q$ by $-q$ in \eqref{r2eq} so that
\begin{align}\label{r2eqa}
\frac{(q^{2};-q^{3})_{\infty}}
{(-q;-q^{3})_{\infty}}
&=
\frac{1}{1}
\+
\frac{q}{1-q}
\+
\frac{q^{3}}{1+q^{2}}
\+
\frac{q^{5}}{1-q^{3}}
\-
\cds \\
&=\frac{1}{1}
\+
\frac{1}{1/q-1}
\+
\frac{1}{1/q^{2}+1}
\+
\frac{1}{1/q^{3}-1}
\-
\cds \notag \\
&=[0;1,1/q-1,1/q^2+1,1/q^{3}-1,1/q^{4}+1,\cdots] \notag
\end{align}
If $q$ is now replaced by $1/a$, where $a>1$ is an integer, and both sides of
\eqref{r2eqa} are inverted, then \eqref{r2a} follows. If $q$ is replaced by $-1/a$
and the resulting  negatives are removed from the continued fraction,
then \eqref{r2b} follows.
\end{proof}

In \cite{ABSYZ03}, the authors generalize the identity at
\eqref{r2eq} as follows.
\begin{proposition}\label{prop2}
Define
\[
F(c,d,q):=\sum_{n=0}^{\infty}\frac{(-1)^{n}c^{n}q^{n(n+1)/2}}{(q;q)_{n}(cq/d;q)_{n}}
\]
and let $\omega = e^{2\pi i/3}$. Then
\begin{equation}\label{r2gen}
\frac{c\,q\,F(c \omega q, \omega^{2},q)}
{(1-c \omega^{2}q)F(c \omega , \omega^{2},q)}=
\frac{cq}{1+cq}
\-
\frac{c^{2}q^{3}}{1+cq^{2}}
\-
\frac{c^{2}q^{5}}{1+cq^{3}}
\-
\cds .
\end{equation}
\end{proposition}
We have changed the statement of their identity slightly to better suit our purposes
and to avoid conflict with our already existing notation.
We use the identity above to prove the following theorem.
\begin{theorem}\label{t4}
Let $a$ be a positive integer and let $c$ be rational such that $a/c \in \mathbb{Z^{+}}$.
If $a/c>1$, then
\begin{equation}\label{t4a}
\left [0;\overline{\frac{a^{2k-1}}{c}-1,\frac{a^{2k}}{c}+1}\right ]_{k=1}^{\infty}=
\frac{c/a\,F(-c \omega /a, \omega^{2},-1/a)}
{(1+c \omega^{2}/a)F(c \omega , \omega^{2},-1/a)}.
\end{equation}
\begin{equation}\label{t4b}
\left [0;\overline{\frac{a^{2k-1}}{c}+1,\frac{a^{2k}}{c}-1}\right ]_{k=1}^{\infty}=
\frac{c/a\,F(c \omega /a, \omega^{2},-1/a)}
{(1-c \omega^{2}/a)F(-c \omega , \omega^{2},-1/a)}.
\end{equation}
If $a/c>2$, then
\begin{equation}\label{t4c}
\left [0;\frac{a}{c}-2,\overline{1,\frac{a^{k+1}}{c}-3}\right ]_{k=1}^{\infty}=
\frac{c/a\,F(-c \omega /a, \omega^{2},1/a)}
{(1+c \omega^{2}/a)F(-c \omega , \omega^{2},1/a)}.
\end{equation}
\begin{equation}\label{t4d}
\left [0;\overline{1,\frac{a^{k}}{c}-1}\right ]_{k=1}^{\infty}=
1-\frac{c/a\,F(c \omega /a, \omega^{2},1/a)}
{(1-c \omega^{2}/a)F(c \omega , \omega^{2},1/a)}.
\end{equation}
\end{theorem}

\begin{proof}
Replace $q$ by $-q$ in \eqref{r2gen} and cancel a negative sign from both sides to get that
\begin{align}\label{r2gena}
\frac{c\,q\,F(-c \omega q, \omega^{2},-q)}
{(1+c \omega^{2}q)F(c \omega , \omega^{2},-q)}&=
\frac{cq}{1-cq}
\+
\frac{c^{2}q^{3}}{1+cq^{2}}
\+
\frac{c^{2}q^{5}}{1-cq^{3}}
\+
\cds \\
&=\frac{1}{1/(cq)-1}
\+
\frac{1}{1/(cq^{2})+1}
\+
\frac{1}{1/(cq^{3})-1}
\+
\cds \notag \\
&=[0; \overline{1/(cq^{2k-1})-1,1/(cq^{2k})+1}]_{k=1}^{\infty}.\notag
\end{align}

If we let $q=1/a$, where $a>1$, $a \in \mathbb{Z^{+}}$ and
$a/c>1$, $a/c \in \mathbb{Z^{+}}$, then \eqref{t4a} follows. If $c$ is replaced
by $-c$ in \eqref{t4a} and both sides are multiplied by $-1$, then \eqref{t4b}
follows.

 If $c$ is replaced by $-c$ and $a$ is replaced by $-a$ in \eqref{t4a}
and the resulting negatives in the continued fraction are cleared, then \eqref{t4c}
follows.

Equation \eqref{t4d} follows from \eqref{t4a} upon replacing $a$ by $-a$,
removing the resulting negatives from the continued fraction and moving the initial $-1$
that appears at the beginning of the continued fraction to the other side of the equation.
\end{proof}

In proving Proposition \ref{prop2} the authors in make use of the following
property of the function $F(a,b,q)$:
\begin{equation*}
\frac{(1-aq/b)F(a,b,q)}{F(aq,b,q)}=
(1-aq-aq/b)-
\frac{a^{2}q^{3}/b}{  \frac{(1-aq^{2}/b)F(aq,b,q)}{F(aq^{2},b,q)}   }.
\end{equation*}

Upon iteration the identity above produces the following continued fraction:
\begin{multline}\label{gen2}
\frac{(1-aq/b)F(a,b,q)}{F(aq,b,q)}=
(1-aq-aq/b)-
\frac{a^{2}q^{3}/b}{1-aq^{2}-aq^{2}/b}
\-\\
\frac{a^{2}q^{5}/b}{1-aq^{3}-aq^{3}/b}
\-
\frac{a^{2}q^{7}/b}{1-aq^{4}-aq^{4}/b}
\-
\cds .
\end{multline}

The authors in \cite{ABSYZ03} did not consider this continued fraction in its full generality,
so did not prove that the right side converged and equalled the left side, but this follows
easily using an argument similar to the one they used for the special case $a=\omega$,
$b=\omega^{2}$. Namely, the continued fraction at \eqref{gen2} is equivalent
to a continued fraction of the form $c_{0}+K_{n=1}^{\infty}c_{n}/1$, where, for $n \geq 2$,
\[
c_{n} = \frac{a^{2}q^{2n+1}}
{(1-aq^{n}-aq^{n}/b)(1-aq^{n+1}-aq^{n+1}/b)}.
\]
Worpitzky's theorem gives that the continued fraction converges. Secondly,
\[
\frac{(1-aq^{n}/b)F(aq^{n-1},b,q)}{F(aq^{n},b,q)} = 1+O(q^{n})
\]
as $n \to \infty$ and this is sufficient to show that the right side of \eqref{gen2}
converges to the left side.
We now use this continued fraction to prove the following theorem.
\begin{theorem}\label{tin}
Let $n>1$ be an integer and $m$ a positive rational such that $m\,n$ is integral.
If $m\,n>4$, then
\begin{equation}\label{tin1}
[mn-3,\overline{ 1,mn^{k} -4}]_{k=2}^{\infty}=
(mn-1) \frac{F(1/m,1,1/n)}{F(1/(mn),1,1/n)}.
\end{equation}
If $m\,n>2$, then
\begin{equation}\label{tin2}
[mn-2,\overline{ mn^{2k} +2, mn^{2k+1} -2}]_{k=1}^{\infty}=
(mn-1) \frac{F(-1/m,1,-1/n)}{F(1/(mn),1,-1/n)}.
\end{equation}
\begin{equation}\label{tin3}
[mn+1,\overline{ 1, mn^{k} }]_{k=2}^{\infty}=
(mn+1) \frac{F(-1/m,1,1/n)}{F(-1/(mn),1,1/n)}.
\end{equation}
If $m\,n^{2}>2$, then
\begin{equation}\label{tin4}
[mn+2,\overline{ mn^{2k} -2, mn^{2k+1} +2}]_{k=1}^{\infty}=
(mn-1) \frac{F(1/m,1,-1/n)}{F(-1/(mn),1,-1/n)}.
\end{equation}
\end{theorem}
\begin{proof}
In \eqref{gen2} above, let $b=1$, divide both sides by $a\,q$ and then
apply a series of similarity transformations to bring all terms into the denominator.
This gives that
\begin{align*}
&\frac{(1/aq-1)F(a,1,q)}{F(aq,1,q)}=
\frac{1}{aq}-2-
\frac{1}{1/aq^{2}-2}
\-
\frac{1}{1/aq^{3}-2}
\-
\frac{1}{1/aq^{4}-2}
\-
\cds \\
&=\frac{1}{aq}-2+
\frac{1}{-1/aq^{2}+2}
\+
\frac{1}{1/aq^{3}-2}
\+
\frac{1}{-1/aq^{4}+2}
\+
\frac{1}{1/aq^{5}-2}
\+
\cds .
\end{align*}
Now let $a=1/n$ and $q=1/n$ and we have that
\begin{equation}\label{baseeq}
(mn-1) \frac{F(1/m,1,1/n)}{F(1/(mn),1,1/n)}
=[mn-2,\overline{ -mn^{2k} +2, mn^{2k+1} -2}]_{k=1}^{\infty}.
\end{equation}
Upon removing the negatives from the left side of \eqref{baseeq}
we get \eqref{tin1}.

Equation \eqref{tin2} follows from \eqref{baseeq}
upon replacing $m$ by $-m$, $n$ by $-n$ and likewise removing the negatives.

Equation \eqref{tin3} follows similarly upon replacing  $n$ by $-n$ and
\eqref{tin4} follows  upon replacing  $m$ by $-m$ (and in each case removing
the negatives).

Remark: The only other value of $b$ which leads to regular continued fraction
expansions is $b=-1$. However, this produce a variant of the generalized
Rogers-Ramanujan continued fraction and, in light of Theorem\ref{t1},
 leads to no new types of regular continued fraction expansions.

\end{proof}

On page 374 of his second notebook, Ramanujan recorded the following continued fraction
identity (see \cite{B98}, page 45, Entry 17).
\begin{proposition}
Let $c$, $d$ and $q$ be complex numbers with $|q|<1$. Define
\[
\phi(c,d,q)= \sum_{n=0}^{\infty}\frac{q^{(n^{2}+n)/2}c^{n}}
{(q;q)_{n}(-dq:q)_{n}}.
\]
Then
\begin{equation}\label{eqt5}
\frac{\phi(c,d,q)}
{\phi(c q,d,q)}=
1+
\frac{cq}{1}
\+
\frac{dq}{1}
\+
\frac{cq^{2}}{1}
\+
\frac{dq^{2}}{1}
\+
\frac{cq^{3}}{1}
\+
\frac{dq^{3}}{1}
\+
\cds .
\end{equation}
\end{proposition}
Once again we have changed the statement of this identity slightly to better suit
our purposes. Equation \ref{eqt5} can be rewritten as
\begin{align}\label{eqt55}
\frac{\phi(c,d,q)}
{\phi(c q,d,q)}&=
1+
\frac{1}{1/(cq)}
\+
\frac{1}{c/d}
\+
\frac{1}{d/(cq)^{2}}
\+
\frac{1}{(c/d)^{2}}
\+
\frac{1}{d^{2}/(cq)^{3}}
\+
\cds \\
&=[1,\overline{d^{k-1}/(cq)^{k},(c/d)^{k}}]_{k=1}^{\infty}. \notag
\end{align}
If we let $c=d m$ and $q=1/(d m n)$, this identity becomes
\begin{equation}\label{eqt55a}
\frac{\phi(d m,d,1/(d m n))}
{\phi(1/n,d,1/(d m n))}=[1,\overline{d^{k-1}n^{k},m^{k}}]_{k=1}^{\infty}.
\end{equation}
From this identity we can deduce the following theorem.
\begin{theorem}\label{t5}
For $r$, $s$ and $q \in \mathbb{C}$ with $|q|<1$, define
\[
\phi(r,s,q)= \sum_{n=0}^{\infty}\frac{q^{(n^{2}+n)/2}r^{n}}
{(q;q)_{n}(-sq;q)_{n}}.
\]
Let $m$ and $n$ be positive integers and let $d$ be rational such that
$d n \in \mathbb{Z^{+}}$ and $d m n>1$. Then
\begin{equation}\label{eqt5a}
[1,\overline{d^{k-1}n^{k},m^{k}}]_{k=1}^{\infty}=
\frac{\phi(d m,d,1/(d m n))}
{\phi(1/n,d,1/(d m n))}.
\end{equation}
If $n>2$ and $m>1$ then
\begin{multline}\label{eqt5b}
[0,\overline{1,d^{2k-2}n^{2k-1}-2,1,m^{2k-1}-1,
d^{2k-1}n^{2k},m^{2k}-1
}]_{k=1}^{\infty}\\
=\frac{\phi(d m,d,-1/(d m n))}
{\phi(-1/n,d,-1/(d m n))}.
\end{multline}
If $m>2$ and $n>1$ then
\begin{multline}\label{eqt5c}
[1,\overline{d^{2k-2}n^{2k-1}-1,1,m^{2k-1}-2,1,
d^{2k-1}n^{2k}-1,m^{2k}
}]_{k=1}^{\infty}\\
=
\frac{\phi(-d m,d,-1/(d m n))}
{\phi(1/n,d,-1/(d m n))}.
\end{multline}
If $n>2$ and $m>1$ then
\begin{multline}\label{eqt5d}
[1,\overline{d^{2k-2}n^{2k-1},m^{2k-1}-1,1,
d^{2k-1}n^{2k}-2,1,m^{2k}-1
}]_{k=1}^{\infty}\\
=
\frac{\phi(-d m,-d,-1/(d m n))}
{\phi(1/n,-d,-1/(d m n))}.
\end{multline}
\end{theorem}
\begin{proof}
Equation \ref{eqt5a} is a restatement of \eqref{eqt55a}, with the stated
conditions on $m$, $n$ and $d$. Equations \ref{eqt5b}, \ref{eqt5c} and
\ref{eqt5d} follow from \ref{eqt5a} by replacing, respectively, $n$ by $-n$,
$m$ by $-m$, $d$ by $-d$ and clearing the negatives from the continued
fraction expansions. Other combinations of sign changes will
produce different regular continued fraction expansions.
\end{proof}

If we replace $q$ by $q^{2}$ and $c$ by $c/q$ in Equation \ref{eqt5} and then apply
Corollary \ref{corcf7}, with once again $c_{1}$ arbitrary, we get that
{\allowdisplaybreaks
\begin{align}\label{eqt555}
\frac{\phi(c/q,d,q^{2})}
{\phi(c q,d,q^{2})}&=
1+
\frac{cq}{1}
\+
\frac{dq^{2}}{1}
\+
\frac{cq^{3}}{1}
\+
\frac{dq^{4}}{1}
\+
\frac{cq^{5}}{1}
\+
\frac{dq^{6}}{1}
\+
\cds\\
=
1-c_{1}+
\frac{c_{1}}{1}
&\-
\frac{cq/c_{1}}{1}
\+
\frac{cq/c_{1}}{1}
\-
\frac{c_{1}dq/c}{1}
\+
\frac{c_{1}dq/c}{1} \notag\\
\-
\frac{(cq)^{2}/(c_{1}d)}{1}
&\+
\frac{(cq)^{2}/(c_{1}d)}{1}
\-
\frac{c_{1}(dq/c)^{2}}{1}
\+
\frac{c_{1}(dq/c)^{2}}{1} \-  \cds \notag\\
\-
\frac{(cq)^{n}/(c_{1}d^{n-1})}{1}
&\+
\frac{(cq)^{n}/(c_{1}d^{n-1})}{1}
\-
\frac{c_{1}(dq/c)^{n}}{1}
\+
\frac{c_{1}(dq/c)^{n}}{1}
\-
\cds \notag\\
=1-c_{1}
\+
\frac{1}{1/c_{1}}
&\-
\frac{1}{c_{1}^{2}/(c q)}
\+
\frac{1}{1/c_{1}}
\-
\frac{1}{c/(d q)} \notag \\
\+
\frac{1}{1/c_{1}}
&\-
\frac{1}{c_{1}^{2}d/(c q)^{2}}
\+
\frac{1}{1/c_{1}}
\-
\frac{1}{(c/(d q))^{2}} \+ \notag \\
&\phantom{assadadsadaasasa}\cds \notag\\
\+
\frac{1}{1/c_{1}}
&\-
\frac{1}{c_{1}^{2}d^{n-1}/(c q)^{n}}
\+
\frac{1}{1/c_{1}}
\-
\frac{1}{(c/(d q))^{n}} \+ \cds .\notag
\end{align}
}
The second equation is valid once again by Worpitzky's Theorem, which
gives that the second continued fraction  converges,
and hence must equal its odd part.

If we now set $c_{1}=1/p$, $c=d\sqrt{m/n}$ and $q=1/\sqrt{mn}$, we have,
once again after a small manipulation to bring all the negative signs into
the denominators, that
\begin{equation}\label{eqt554}
\frac{\phi(d m,d,1/(m n))}
{\phi(d/n,d,1/(m n))}
=-\frac{1}{p}+
\left [1; \overline{p, \frac{-n^{k}}{p^{2}d},-p, m^{k}} \, \right ]_{k=1}^{\infty}.
\end{equation}
As with previous identities, we now choose various combinations of
signs for the variables to get the following theorem.
\begin{theorem}\label{t6}
Let $d$ be a positive rational and let $m$, $n$ and $p$ be integers such that $m n >1$
and $p^{2}d|n$. If $m>1$, $n/(p^{2}d)>1$ and $p>1$, then
\begin{equation}\label{t6a}
\left [1; \overline{p-1,1, \frac{n^{k}}{p^{2}d}-1,p-1,1, m^{k}-1} \, \right ]_{k=1}^{\infty}
=\frac{1}{p}+\frac{\phi(d m,d,1/(m n))}
{\phi(d/n,d,1/(m n))}.
\end{equation}
If $m>1$, $n/(p^{2}d)>1$ and $p>2$, then
\begin{equation}\label{t6b}
\left [1; \overline{p, \frac{n^{k}}{p^{2}d}-1,1,p-2,1, m^{k}-1} \, \right ]_{k=1}^{\infty}
=\frac{1}{p}+\frac{\phi(-d m,-d,1/(m n))}
{\phi(-d/n,-d,1/(m n))}.
\end{equation}
If $m>1$, $n/(p^{2}d)>1$ and $p>1$, then
\begin{equation}\label{t6c}
\left [0; \overline{1,p-1, \frac{n^{k}}{p^{2}d}-1,1,p-1, m^{k}-1} \, \right ]_{k=1}^{\infty}
=\frac{-1}{p}+\frac{\phi(d m,d,1/(m n))}
{\phi(d/n,d,1/(m n))}.
\end{equation}
If $m>1$, $n/(p^{2}d)>1$ and $p>2$, then
\begin{multline}\label{t6d}
\text{\small{$
\left [1; \overline{p-1,1, \frac{n^{2k-1}}{p^{2}d}-1,p, m^{2k-1}-1,1,
p-2,1, \frac{n^{2k}}{p^{2}d}-1,p-1,1, m^{2k}-1} \, \right ]_{k=1}^{\infty}$}}\\
=\frac{1}{p}+\frac{\phi(-d m,d,-1/(m n))}
{\phi(d/n,d,-1/(m n))}.
\end{multline}
If $m>1$, $n/(p^{2}d)>1$ and $p>2$, then
\begin{multline}\label{t6e}
\text{\small{$
\left [1; \overline{p,\frac{n^{2k-1}}{p^{2}d}-1,1,p-2,1, m^{2k-1}-1,
p-1,1, \frac{n^{2k}}{p^{2}d}-1,p-1,1, m^{2k}-1} \, \right ]_{k=1}^{\infty}$}}\\
=\frac{1}{p}+\frac{\phi(d m,d,-1/(m n))}
{\phi(-d/n,d,-1/(m n))}.
\end{multline}
\end{theorem}

\begin{proof}
Equations \ref{t6a} to \ref{t6e} follow from Equation \ref{eqt554}
by, respectively, i) keeping the sign of all variables unchanged,
ii) replacing $d$ by $-d$, iii) replacing $p$ by $-p$,
iv) replacing $m$ by $-m$, v) replacing $n$ by $-n$ and in each case
removing the negative signs from the resulting continued fraction.

Other combinations of signs will lead to other regular continued
 fraction expansions.
\end{proof}

\section{Concluding Remarks}
Several of the continued fraction identities stated by Ramanujan
have been generalized by various authors and it is possible that some of these
can also be manipulated to give new classes of regular continued fraction
expansions whose limits can represented in other ways.

In \cite{DNNS97} the authors prove that if $0<|q|<1$ and is algebraic, then the
 Rogers-Ramanujan   continued fraction $K(q)$ at \eqref{rreq} converges
to a transcendental number. This means that, in Theorem \ref{t1},
if we set $c=a^j$ and $d=a^{k}$,
where $j$ and $k$ are non-negative integers,
  then the values of all of the resulting continued fractions are transcendental.
Similarly,  setting $r=s$ in Theorem \ref{t2} gives that the
values of all of the resulting continued fractions are
transcendental.   It is likely that transcendence results for
other $q$-continued fractions may be used to show the
transcendence of the values of some of the other regular continued
fractions in the paper, but we have not pursued that here.

\allowdisplaybreaks{
}

\end{document}